\documentclass[a4paper, 10pt]{article}
\usepackage{amssymb, amsmath, amsfonts}
\usepackage[english]{babel}
\usepackage{amsthm,amscd,amsfonts,latexsym,color,epsfig}

\begin{document}

\begin{center}
\sloppy
\textbf{\Large Necessary non-local conditions for a time-fractional diffusion-wave equation}
\bigskip

\centerline {{\bf Murat Mamchuev}
\footnote{
Institute of Applied Mathematics and Automation of KBSC of RAS, 
Shortanova str. 89-A, Nalchik, 360000, Kabardino-Balkar Republic, Russia, E-mail: mamchuev@rambler.ru}}
\end{center}

\bigskip

\begin{abstract}
In this paper the time-fractional diffusion-wave equation with Riem\-man-Liouville  fractional deri\-va\-tive is studied.
The integral operators with the Wright function in the kernel, associated with the studied equation are introdused and their properties are investigated.
In terms of these operators the necessary non-local conditions binding traces of solution and its derivatives on the boundary of a rectangular domain are found. 
By using the limiting properties of the Wright function the necessary non-local conditions  for wave equation are obtained.
With the help of the mentioned integral operator's properties a unique solvability of the problem with Samarskii integral condition 
for the diffusion-wave equation is proved. 
The solution is obtained in explicit form.
\end{abstract}

\vspace{5mm}

\section{Introduction}

In the domain $\Omega=\{(x,y):\ 0<x<l,\ 0<y<T\}$ we consider the equation
\begin{equation} \label{lmamud}
u_{xx}(x,y)-D_{0y}^{\alpha} u(x,y) =0,
\end{equation}
where
$0<\alpha<2,$
$D_{ay}^{\nu}$ is the Riemann-Liouville fractional integro-differentia\-tion operator
of order $\nu$  \cite [p. 9]{n}, which determined as 
$$D_{ay}^{\nu}g(y)=\frac{\mathop{\rm sgn}(y-a)}{\Gamma(-\nu)}
\int\limits_a^y \frac{g(s)ds}{|y-s|^{\nu+1}}, \quad  \nu< 0,$$
for $\nu \geq 0$ the operator $D_{ay}^{\nu}$ can be determined by
recursive relation
$$D_{ay}^{\nu}g(y)=\mathop{\rm sgn}(y-a)\frac{d}{dy}D_{ay}^{\nu-1}g(y), \quad \nu \geq 0,
$$
$\Gamma(z)$ is the gamma-function.

It is well-known that the fractional partial differential equations
are appire in mathematical
models describing various processes in a mediums with a fractal structure  
(see for example \cite [chapter 5]{n}).
Equations of the form (\ref{lmamud}), due to their numerous applications in the modeling of processes of a different nature, 
have been actively investigated by many authors using different methods during  the last decades.
Such equations describe anomalous diffusion and subdiffusion processes, 
relaxation phenomena in complex viscoelastic materials, and so on.

For an extensive  bibliography on this subject see for example references in papers 
\cite{PsMon},  \cite{pizvran},
\cite{MamMon}, \cite{Mamchuev-2017}.

\section{Integral operators associated with \\ the diffusion-wave equation }

In the paper \cite{g1} (see also \cite{Mainardi-1996}) the fundamental solution $\Gamma(x-t,y-s)$ 
of the equation (\ref{lmamud}) was represented in term of the function
\begin{equation} \label{nnfr0}
\Gamma(x,y)
=\frac{y^{\beta-1}}{2}
\phi(-\beta,\beta;-|x|y^{-\beta}),
\end{equation}
where 
$\beta=\alpha/2$, $\phi(\rho,\mu;z)=\sum\limits_{k=0}^{\infty}\frac{z^k}{k!\Gamma(\rho k+\mu)}$ is the Wright function \cite{Wright-1933}.

Let us introduse the following operators ${\mathcal N}_{x_1x_2}^{\theta,x,y}$ и
${\mathcal R}_{0y}^{\delta,x},$ which acting by the formulas
\begin{equation} \label{lmamna}
{\mathcal N}_{x_1x_2}^{\theta,x,y}\nu(t)=\frac{1}{2}\int\limits_{x_1}^{x_2}
\nu(t)y^{\theta-1}
\phi(-\beta,\theta;-|x-t|y^{-\beta})dt,
\end{equation}
\begin{equation} \label{lmamrn}
{\mathcal R}_{0y}^{\delta,x}\mu(y)=\frac{1}{2}\int\limits_{0}^{y}
\mu(s)(y-s)^{\delta-1}
\phi(-\beta,\delta;-|x|(y-s)^{-\beta})ds.
\end{equation}

The following properties of the operators ${\mathcal N}_{x_1x_2}^{\theta,x,y}$ and
${\mathcal R}_{0y}^{\delta,x}$ are hold.

{\bf Property 1.} { The relation}
\begin{equation} \label{ls1}
\left(2{\mathcal R}_{0y}^{\delta,x_1}\right)
\left(2{\mathcal R}_{0y}^{\theta,x_2}\right)\mu(y)=
2{\mathcal R}_{0y}^{\delta+\theta,x_1+x_2}\mu(y)
\end{equation}
holds for all $x_1>0,$ $x_2>0.$

{\bf Property 2.} { Let $a\geq 0,$ $\delta+\beta>0,$ $\theta+\beta>0.$
Then the following relations hold}
$$
{\mathcal R}_{0y}^{\delta,a}{\mathcal N}_{0l}^{\theta,b,y}\tau(t)=
\frac{1}{2}\left\{
\begin{array}{ll}
{\mathcal N}_{0l}^{\delta+\theta,b-a,y}\tau(t),    & b\leq 0,\\
\left[{\mathcal N}_{0b}^{\delta+\theta,a+b,y}+
{\mathcal N}_{bl}^{\delta+\theta,b-a,y}\right]\tau(t),    & 0<b<l, \\
{\mathcal N}_{0l}^{\delta+\theta,a+b,y}\tau(t),    & b\geq l. \\
\end{array}
\right.
$$

{\bf Property 3.}
{ Let $y^{1-\nu}\mu(y)\in C[0,T],$ $\nu+\delta\geq 0,$ then }
\begin{equation} \label{ls6}
\lim\limits_{x\rightarrow 0}{\mathcal R}_{0y}^{\delta,x} \mu(y)=
\frac{1}{2}D_{0y}^{-\delta}\mu(y).
\end{equation}

{\bf Property 4.} {Let $\tau(t)\in C[0,T],$ then }
\begin{equation} \label{ls5}
\lim\limits_{y\rightarrow 0}D_{0y}^{2\beta-n}
{\mathcal N}_{0l}^{\beta-k+1,x,y}\tau(x)=
\left\{
\begin{array}{ll}
\tau(x),& k=n, \\
0,	& k<n.
\end{array}
\right.
\end{equation}

{\bf Property 5.} { The relation}
$$\int\limits_0^l{\mathcal R}_{0y}^{\theta,a\pm x} \mu(y)dx=
\pm\left[{\mathcal R}_{0y}^{\theta+\beta,a}-
{\mathcal R}_{0y}^{\theta+\beta,a\pm l}\right] \mu(y)$$
holds.

{\bf Property 6.} { The folloving relations hold}
$$\int\limits_0^l{\mathcal N}_{0l}^{\delta,a\pm x,y} \tau(t)dx=
\pm{\rm sign\,} a\left[{\mathcal N}_{0l}^{\delta+\beta,a,y}-
{\mathcal N}_{0l}^{\delta+\beta,a\pm l,y}\right] \tau(t), \quad |a|\geq 2l,$$
$$\int\limits_0^l{\mathcal N}_{0l}^{\delta,-x,y} \tau(t)dx=
\left[{\mathcal N}_{0l}^{\delta+\beta,0,y}-
{\mathcal N}_{0l}^{\delta+\beta,-l,y}\right] \tau(t),
$$
$$\int\limits_0^l{\mathcal N}_{0l}^{\delta,x,y} \tau(t)dx=
-\left[{\mathcal N}_{0l}^{\delta+\beta,0,y}+
{\mathcal N}_{0l}^{\delta+\beta,l,y}\right] \tau(t)+
\frac{y^{\delta+\beta-1}}{\Gamma(\delta+\beta)}\int\limits_0^l\tau(x)dx.
$$

For the proof of properties  1 and 2 it is enough to change the integration order, use the 
deffenitions of the operators
${\mathcal R}_{0y}^{\delta,x}$ и ${\mathcal N}_{0l}^{\theta,x,y}$
and 
the following convolution formula for Wright functions \cite{Stankovic-1970}
$$\int\limits_{0}^{y}
\xi^{\delta-1}\phi(-\alpha,\delta;-x_1 \xi^{-\alpha})
(y-\xi)^{\mu-1}\phi(-\alpha,\mu;-x_2(y-\xi)^{-\alpha})d\xi=
$$
\begin{equation}\label{veq10}
=y^{\delta+\mu-1}\phi(-\alpha,\delta+\mu;-(x_1+x_2) y^{-\alpha}),
\quad \forall x_1, x_2>0.
\end{equation}

The property 3 is proved in \cite[p. 35]{PsMon}.
We rewrited it in the terms of the operator ${\mathcal R}_{0y}^{\delta,x}.$

The property 4 follows from the formula  \cite[p. 24]{PsMon}
\begin{equation} \label{vnn}
D^\nu_{0y}y^{\delta-1}\phi(-\beta,\delta;-cy^{-\beta})=
y^{\delta-\nu-1}\phi(-\beta,\delta-\nu;-cy^{-\beta}),
\quad c>0,
\quad \delta+\beta>0,
\end{equation}
the relation \cite{Mainardi-1996}
\begin{equation}\label{vf15}
\int\limits_{0}^{\infty}\phi(-\beta,1-\beta;-s)ds=1,
\end{equation}
and the estimates \cite[p. 29]{PsMon}
\begin{equation} \label{vofwr0}
\left|y^{\delta-1}\phi(-\beta,\delta;-x y^{-\beta})\right|
\leq C x^{-\theta}y^{\delta+\beta\theta-1}, \quad x>0, \quad y>0,
\end{equation}
where $\theta\geq 0$ for $\delta\notin {\mathbb N}_0;$ and
$\theta\geq -1$ for $\delta\in {\mathbb N}_0;$ ${\mathbb N}_0={\mathbb N}\cup \{0\},$
$C$ is a positive constant. 

For the proof of properties 5 and 6 it is enough to use the deffinitions of operators
${\mathcal R}_{0y}^{\delta,x}$ and ${\mathcal N}_{0l}^{\theta,x,y},$
and formula \cite{Wright-1933}
\begin{equation}\label{vf13}
\frac{d}{dz}\phi(-\beta, \delta; z)=\phi(-\beta, \delta-\beta; z).
\end{equation}

In view of (\ref{ls6}) from properties 1 and 2 follow the relations
\begin{equation} \label{ls2}
D_{0y}^{-\theta}{\mathcal R}_{0y}^{\delta,x}\mu(y)=
{\mathcal R}_{0y}^{\delta,x}D_{0y}^{-\theta}\mu(y)
={\mathcal R}_{0y}^{\delta+\theta,x}\mu(y),
\end{equation}
\begin{equation} \label{ls3}
D_{0y}^{-\delta}{\mathcal N}_{0l}^{\theta,x,y}\nu(x)=
{\mathcal N}_{0l}^{\delta+\theta,x,y}\nu(x).
\end{equation}

\section{Necessary non-local conditions for the \\ time-fractional diffusion-wave equation}

Let  $n\in \{1,2\}$ is a number such that $n-1<\alpha\leq n.$
A function $u=u(x,y)$ of the class $D_{0y}^{\alpha-k}u(x,y)\in C(\bar{\Omega}),$ $1\leq k\leq n,$
$u_{xx}(x,y),$ $D_{0y}^{\alpha}u(x,y)\in C(\Omega),$ satisfying the equation (\ref{lmamud}) at all points $(x,y)\in \Omega$ 
is called {\it a regular solution of equation} (\ref{lmamud}) in the domain $\Omega$ 
\cite [p. 103]{PsMon}.

The following assertion hold \cite{Mamchuev-2007}.

{\bf Theorem 1.} { Let $u=u(x,y)$ is a regular in the domain $\Omega$ solution of equation (\ref{lmamud}),
satisfying the condition 
\begin{equation} \label{lmamnu}
\lim\limits_{y\rightarrow 0} D_{0y}^{\alpha-k}u(x,y)=\tau_k(x),
\quad 1\leq k\leq n,
\quad 0<x<l,
\end{equation}
such that $u_x\in C([0,l]\times (0,T))$
and $u_x(0,y),\, u_x(l,y)\in L[0,T].$
Then the function $u(x,y)$ fulfill the non-local conditions
\begin{equation} \label{lmamx0}
u(0,y)=2\sum\limits_{k=1}^{n}{\mathcal N}_{0l}^{\beta-k+1,0,y}\tau_k(t)+
2{\mathcal R}_{0y}^{\beta,l}u_x(l,s)-
D_{0y}^{-\beta}u_x(0,s)
+2{\mathcal R}_{0y}^{0,l}u(l,s),
\end{equation}
\begin{equation} \label{lmamxl}
u(l,y)=2\sum\limits_{k=1}^{n}{\mathcal N}_{0l}^{\beta-k+1,l,y}\tau_k(t)-
2{\mathcal R}_{0y}^{\beta,l}u_x(0,s)
+D_{0y}^{-\beta}u_x(l,s)
+2{\mathcal R}_{0y}^{0,l}u(0,s).
\end{equation}
}

{\bf Proof.}
We use the general representation of regular solutions of equation (\ref{lmamud})
in the rectangular domain \cite[c.116]{PsMon}
$$u(x,y)=\sum\limits_{k=1}^{n}(-1)^{k-1}\int\limits_{x_1}^{x_2}
\tau_k(t)\frac{\partial^{k-1}}{\partial s^{k-1}}G(x,y;t,0)dt+$$
\begin{equation} \label{lmamfg1}
+\sum\limits_{i=1}^{2}(-1)^{i}\int\limits_{0}^{y}
[G(x,y;x_i,s)u_{t}(x_i,s)-G_{t}(x,y;x_i,s)u(x_i,s)]ds.
\end{equation}
By setting in relation (\ref{lmamfg1})
$G(x,y;t,s)=\Gamma(x-t,y-s),$
taking into account the equalities
$$\frac{\partial}{\partial t}\Gamma(x,y)=\frac{{\rm sgn} \,x}{2y}\phi(-\beta,0;-|x|y^{-\beta}),
$$
$$\frac{\partial}{\partial s}\Gamma(x,y)=-\frac{y^{\beta-2}}{2}\phi(-\beta,\beta-1;-|x|y^{-\beta}),
$$
which hold in view of relations (\ref{vnn}) and (\ref{vf13}),
we obtain
$$u(x,y)=\frac{1}{2}\sum\limits_{k=1}^{n}
\int\limits_{x_1}^{x_2}\tau_k(t)y^{\beta-k}
\phi(-\beta,\beta-k+1;-|x-t|y^{-\beta})dt+$$
$$+\frac{1}{2}\int\limits_{0}^{y}u_{t}(x_2,s)(y\!-\!s)^{\beta-1}
\phi(-\beta,\beta;-(x_2-x)(y-s)^{-\beta})ds-$$
$$-\frac{1}{2}\int\limits_{0}^{y}u_{t}(x_1,s)(y\!-\!s)^{\beta-1}
\phi(-\beta,\beta;-(x-x_1)(y-s)^{-\beta})ds+$$
$$+\frac{1}{2}\int\limits_{0}^{y}\frac{u(x_2,s)}{y-s}
\phi(-\beta,0;-(x_2-x)(y-s)^{-\beta})ds+
$$
\begin{equation} \label{lmamfg2}
+\frac{1}{2}\int\limits_{0}^{y}\frac{u(x_1,s)}{y-s}
\phi(-\beta,0;-(x-x_1)(y-s)^{-\beta})ds.
\end{equation}
In the terms of the operators (\ref{lmamna}) and (\ref{lmamrn}), the relation (\ref{lmamfg2}) can be rewriten in the form
$$u(x,y)=\sum\limits_{k=1}^{n}{\mathcal N}_{x_1x_2}^{\beta-k+1,x,y}\tau_k(t)+
{\mathcal R}_{0 y}^{\beta,x_2-x}u_{t}(x_2,s)-
{\mathcal R}_{0 y}^{\beta,x-x_1}u_{t}(x_1,s)+$$
\begin{equation} \label{lmamfg3}
+{\mathcal R}_{0 y}^{0,x_2-x}u(x_2,s)
+{\mathcal R}_{0 y}^{0,x-x_1}u(x_1,s).
\end{equation}
By using the fact that in view of the property 3 of the operator
${\mathcal R}_{0y}^{\delta,x}$
$$
\lim\limits_{x\rightarrow 0}{\mathcal R}_{0y}^{\beta,x}\mu(y)=\frac{1}{2}D_{0y}^{-\beta}\mu(y),
\quad
\lim\limits_{x\rightarrow 0}{\mathcal R}_{0y}^{0,x}\mu(y)=\frac{\mu(y)}{2},
$$
from (\ref{lmamfg3}) we obtain
$$\frac{1}{2}u(x_1,y)=
\sum\limits_{k=1}^{n}{\mathcal N}_{x_1x_2}^{\beta-k+1,x_1,y}\tau_k(t)+
{\mathcal R}_{0y}^{\beta,x_2-x_1}u_{t}(x_2,s)-
$$
\begin{equation} \label{lmamx1}
-\frac{1}{2}D_{0y}^{-\beta}u_{t}(x_1,s)+
{\mathcal R}_{0y}^{0,x_2-x_1}u(x_2,s),
\end{equation}
$$\frac{1}{2}u(x_2,y)=\sum\limits_{k=1}^{n}{\mathcal N}_{x_1x_2}^{\beta-k+1,x_2,y}\tau_k(t)-
{\mathcal R}_{0y}^{\beta,x_2-x_1}u_{t}(x_1,s)+
$$
\begin{equation} \label{lmamx2}
+\frac{1}{2}D_{0y}^{-\beta}u_{t}(x_2,s)+
{\mathcal R}_{0y}^{0,x_2-x_1}u(x_1,s).
\end{equation}
By passing in the relations (\ref{lmamx1}) and (\ref{lmamx2}) to the limit as 
$(x_1,x_2)\rightarrow (0,l),$ we obtain  (\ref{lmamx0}) and (\ref{lmamxl}).
The proof of Theorem 1 is complete.

For $\alpha=1$ the conditions (\ref{lmamx0}) and (\ref{lmamxl}) are coinside with necessary
non-local conditions for the heat equation \cite [p. 275]{n1}.

\section{Necessary non-local conditions \\ for the wave equation}

The following lemma holds \cite{pizvran}.

{\bf Lemma 1.} { Let the function $g(t)$ is absolutly integrable 
on any finite interval of the semiaxis
$t>0,$  is continuous at the point $t=1$ and grows with $t\to \infty$ no faster than
${\rm exp}\{\sigma t^{\delta}\},$ $\sigma>0,$
$\delta<\frac{1}{1-\beta}.$
Then }
$$\lim\limits_{\beta\to 1}\int\limits_{0}^{\infty}
g(t)\phi(-\beta,0;-t)dt=g(1), \quad
\lim\limits_{\beta\to 1}\int\limits_{0}^{\infty}
g(t)\phi(-\beta,\beta;-t)dt=\int\limits_{0}^{1}g(t)dt.$$

Using Lemma 1, we proceed in the equalities (\ref{lmamx0}) и (\ref{lmamxl}) to the limit for $\alpha\to 2.$
For this purpose we rewrite conditions (\ref{lmamx0}) and (\ref{lmamxl}) in the following form
$$u(0,y)=\int\limits_0^l\tau_1(\xi)y^{\beta-1}
\psi(-\beta,\beta;-\xi y^{-\beta})d\xi+
\int\limits_0^l\tau_2(\xi)y^{\beta-2}
\phi(-\beta,\beta-1;-\xi y^{-\beta})d\xi+$$
$$+\int\limits_0^y\frac{u_x(l,\eta)}{(y-\eta)^{1-\beta}}
\phi(-\beta,\beta;-l (y-\eta)^{-\beta})d\eta-
\frac{1}{\Gamma(\beta)}\int\limits_0^y\frac{u_x(0,\eta)}{(y-\eta)^{1-\beta}}d\eta+
$$
\begin{equation} \label{lu0y}
+\int\limits_0^y\frac{u(l,\eta)}{y-\eta}
\phi(-\beta,0;-l(y-\eta)^{-\beta})d\eta
=\sum\limits_{i=1}^{5}I_i(y),
\end{equation}
$$u(l,y)=\int\limits_0^l\tau_1(\xi)y^{\beta-1}
\phi(-\beta,\beta;-(l-\xi)y^{-\beta})d\xi+$$
$$+\int\limits_0^l\tau_2(\xi)y^{\beta-2}
\phi(-\beta,\beta-1;-(l-\xi)y^{-\beta})d\xi-$$
$$-\int\limits_0^y\frac{u_x(0,\eta)}{(y-\eta)^{1-\beta}}
\phi(-\beta,\beta;-l(y-\eta)^{-\beta})d\eta+
\frac{1}{\Gamma(\beta)}\int\limits_0^y\frac{u_x(l,\eta)}{(y-\eta)^{1-\beta}}d\eta+
$$
\begin{equation} \label{luly}
+\int\limits_0^y\frac{u(0,\eta)}{y-\eta}
\phi(-\beta,0;-l(y-\eta)^{-\beta})d\eta
=\sum\limits_{i=1}^{5}J_i(y).
\end{equation}
Next, we transform the integrals $I_i(y)$
$$I_1(y)=\int\limits_0^l\tau_1(\xi)y^{\beta-1}
\phi(-\beta,\beta;-\xi y^{-\beta})d\xi=
\int\limits_0^{l/y^{\beta}}\tau_1(y^{\beta}\eta)y^{2\beta-1}
\phi(-\beta,\beta;-\eta)d\eta=$$
$$=\int\limits_0^{\infty}\tau_1(y^{\beta}\eta)y^{2\beta-1}
\phi(-\beta,\beta;-\eta)H(l-y^{\beta}\eta)d\eta,
$$
$$I_2(y)=\int\limits_0^l\tau_2(\xi)y^{\beta-2}
\phi(-\beta,\beta-1;-\xi y^{-\beta})d\xi=
\frac{d}{dy}\int\limits_0^l\tau_2(\xi)y^{\beta-1}
\phi(-\beta,\beta;-\xi y^{-\beta})d\xi,
$$
$$I_3(y)\!=\!\!\int\limits_0^y\!\frac{u_x(l,\eta)}{(y-\eta)^{1-\beta}}
\phi(-\beta,\beta;-l(y-\eta)^{-\beta})d\eta=$$
$$=\!\int\limits_{l/y^{\beta}}^{\infty}\!u_x\left(l,y-(l/\xi)^{\frac{1}{\beta}}\right)
\frac{l}{\beta\xi^2}\phi(-\beta,\beta;-\xi)d\xi=$$
$$=\int\limits_0^{\infty}u_x\left(l,y-(l/\xi)^{\frac{1}{\beta}}\right)
\frac{l}{\beta\xi^2}\phi(-\beta,\beta;-\xi)H(\xi-l/y^{\beta})d\xi,
$$
$$I_5(y)=\int\limits_0^y\frac{u(0,\eta)}{y-\eta}
\phi(-\beta,0;-l(y-\eta)^{-\beta})d\eta =$$
$$
=\int\limits_{0}^{\infty}u\left(l,y-(l/\xi)^{\frac{1}{\beta}}\right)
\frac{1}{\beta\xi}\phi(-\beta,0;-\xi)H(\xi-l/y^{\beta})d\xi.$$
Using Lemma 1, we obtain the following relations
$$\lim\limits_{\beta\to 1}I_1(y)=\int\limits_{0}^{1}\tau_1(y\eta)H(l-y\eta)yd\eta=
\int\limits_{0}^{y}\tau_1(\xi)H(l-\xi)d\xi=
$$
\begin{equation} \label{li1}
=\int\limits_{0}^{l}\tau_1(\xi)H(y-\xi)d\xi,
\end{equation}
\begin{equation} \label{li2}
\lim\limits_{\beta\to 1}I_2(y)=\frac{d}{dy}
\int\limits_{0}^{y}\tau_2(\xi)H(l-\xi)d\xi=\tau_2(y)H(l-y),
\end{equation}
$$\lim\limits_{\beta\to 1}I_3(y)=
\int\limits_{0}^{1}u_x(l,y-l/\xi)\frac{l}{\xi^2}H(\xi-l/y)d\xi=
$$
\begin{equation} \label{li3}
=H(y-l)\int\limits_{l/y}^{1}u_x(l,y-l/\xi)\frac{l}{\xi^2}d\xi=
H(y-l)\int\limits_{0}^{y-l}u_x(l,\eta)d\eta,
\end{equation}
\begin{equation} \label{li5}
\lim\limits_{\beta\to 1}I_5(y)=
u(l,y-l)H(y-l),
\end{equation}
where  $H(x)$ is the Heaviside function.

Similarly, for the terms on the right-hand side of (\ref{luly}), we obtain
\begin{equation} \label{lj1}
\lim\limits_{\beta\to 1}J_1(y)=\int\limits_{l-y}^{l}\tau_1(\xi)H(\xi)d\xi,
\end{equation}
\begin{equation} \label{lj2}
\lim\limits_{\beta\to 1}J_2(y)=\frac{d}{dy}
\int\limits_{l-y}^{l}\tau_2(\xi)H(\xi)d\xi=\tau_2(l-y)H(l-y),
\end{equation}
\begin{equation} \label{lj3}
\lim\limits_{\beta\to 1}J_3(y)=
H(y-l)\int\limits_{0}^{y-l}u_x(0,\eta)d\eta,
\end{equation}
\begin{equation} \label{lj5}
\lim\limits_{\beta\to 1}J_5(y)=
u(0,y-l)H(y-l).
\end{equation}

Using relations (\ref{li1}) -- (\ref{lj5}), from (\ref{lu0y}) and (\ref{luly})
we obtain necessary non-local condition for wave equation
$$u(0,y)=\int\limits_{0}^{l}\tau_1(t)H(y-t)dt+\tau_2(y)H(l-y)
+H(y-l)\int\limits_{0}^{y-l}u_x(l,s)ds-$$
$$-\int\limits_{0}^{y}u_x(0,s)ds
+u(l,y-l)H(y-l),$$
$$u(l,y)=\int\limits_{l-y}^{l}\tau_1(t)H(t)dt+\tau_2(l-y)H(l-y)
-H(y-l)\int\limits_{0}^{y-l}u_x(0,s)ds+$$
$$+\int\limits_{0}^{y}u_x(l,s)ds
+u(0,y-l)H(y-l).$$
For $y\leq l$ the conditions take the form:
\begin{equation} \label{lmamu11}
u(0,y)=\int\limits_{0}^{y}\tau_1(t)dt+\tau_2(y)-\int\limits_{0}^{y}u_x(0,s)ds,
\end{equation}
\begin{equation} \label{lmamu12}
u(l,y)=\int\limits_{l-y}^{l}\tau_1(t)dt+\tau_2(l-y)+\int\limits_{0}^{y}u_x(l,s)ds.
\end{equation}
For $y\geq l$ we have
\begin{equation} \label{lmamu31}
u(0,y)=\int\limits_{0}^{l}\tau_1(t)dt+\int\limits_{0}^{y-l}u_x(l,s)ds-
\int\limits_{0}^{y}u_x(0,s)ds+u(l,y-l),
\end{equation}
\begin{equation} \label{lmamu32}
u(l,y)=\int\limits_{0}^{l}\tau_1(t)dt-\int\limits_{0}^{y-l}u_x(0,s)ds+
\int\limits_{0}^{y}u_x(l,s)ds+u(0,y-l).
\end{equation}

In case $l=T$ the conditions (\ref{lmamu11}) and (\ref{lmamu12}) was obtained
in paper \cite{zn2006}.

\section{Samarskii problem}

\subsection{Samarskii problem for the time-fractional diffusion-wave equation}

Using Theorem 1, we investigate the Samarskii problem for the equation (\ref{lmamud}) 
in the following formulation.

{\bf Problem 1.}
{ In the domain $\Omega$ find a solution $u(x,y)$ of equation (\ref{lmamud}),
satisfying the condition (\ref{lmamnu}) and boundary conditions 
\begin{equation} \label{pu2}
a_1u(0,y)+a_2u(l,y)=\varphi(y),\  0< y\leq T,
\end{equation}
\begin{equation} \label{pu3}
\int\limits_0^lu(x,y)dx=\mu(y),\   0< y \leq T,
\end{equation}
where $\tau_k(x),$ $\varphi (y),$ $\mu (y)$ are given functions, $a_1,$ $a_2$ are
given numbers, and $a_1\not=a_2.$ }

The condition (\ref{pu3}) is called Samarskii condition \cite[p. 140]{n1}.

In paper \cite{Bazhlekova-1996}  Duhamel-type representation of the solution of problem 1
with fractional derivative in Caputo sence was obtained
in case $a_1=1,$ $a_2=0,$ 
$\varphi\equiv 0,$
$\mu\equiv 0,$
and the initial conditions given in the form
$$\frac{\partial^k}{\partial y^k}u(x,y)\big|_{y=0}=f_k(x), \quad k=1,n,$$
where $f_k(x)$ 
are the functions choosed in a special way. 
In paper \cite {zn1997} in case $0<\alpha <1,$ $a_1=1,$ $a_2=0,$
$\mu(y)\equiv {\mathop{\rm const}}\,y^{\alpha-1},$
and when the condition (\ref{lmamnu}) given in local form,
problem 1 was studied by using the separation of variables method.
Note that the solving of Problem 1 by the reduction to the first boundary value problem 
with the help of conditions (\ref{lmamx0}) and (\ref{lmamxl}) was announced in paper \cite{Mamchuev-2007}.

By $C^{1,q}[0,l]$  we denote the space of continuously
differentiable functions on $[0,l]$ whose derivatives satisfy the Holder condition with exponent $q.$

Following assertion hold.

{\bf Theorem 2.}
{ Let
$\tau_1(x)\in C[0,l];$
$\tau_2(x)\in C^{1,q}[0,l],$  $q>\frac{1-\beta}{\beta},$
for $n=2;$
$y^{n-\alpha}\varphi(y)\in C[0,T],$ $D_{0y}^{\alpha}\mu(y)\in C[0,T]$
and the matching conditions 
\begin{equation} \label{pus1}
\lim\limits_{y\rightarrow 0}D_{0y}^{\alpha-n}\varphi(y)=a_1\tau_n(0)+a_2\tau_n(l),
\end{equation}
\begin{equation} \label{pus2}
\lim\limits_{y\rightarrow 0}D_{0y}^{\alpha-k}\mu(y)=\int\limits_0^l\tau_k(x)dx,
\quad k=1,n.
\end{equation}
are satisfied.
Then there exists unique regular in the domain $\Omega$ solution of problem 1.
This solution has the form
$$u(x,y)=\sum\limits_{k=1}^{n}\sum\limits_{m=-\infty}^{\infty}
\left[{\mathcal N}_{0l}^{\beta-k+1,2ml+x,y}-{\mathcal N}_{0l}^{\beta-k+1,2ml-x,y}\right]
\tau_k(x)-$$
\begin{equation} \label{prz1}
-2\sum\limits_{m=1}^{\infty}
\left[{\mathcal R}_{0y}^{0,2ml-l+x}-
{\mathcal R}_{0y}^{0,2ml-l-x}\right]\varphi_l(y)
+2\sum\limits_{m=1}^{\infty}
\left[{\mathcal R}_{0y}^{0,2ml-2l+x}-{\mathcal R}_{0y}^{0,2ml-x}\right]\varphi_0(y),
\end{equation}
where
$$\varphi_0(y)=\frac{a_2}{a_2-a_1}\psi(y)+\frac{1}{a_2-a_1}\varphi(y),
\quad
\varphi_1(y)=\frac{a_1}{a_1-a_2}\psi(y)+\frac{1}{a_1-a_2}\varphi(y),
$$
\begin{equation} \label{ppsi}
\psi(y)=2\sum\limits_{k=1}^{n}\sum\limits_{m=-\infty}^{\infty}
{\mathcal N}_{0l}^{\beta-k+1,ml,y}\tau_k(\xi)
+4\sum\limits_{m=1}^{\infty}{\mathcal R}_{0y}^{\beta,ml}D_{0y}^{\alpha}\mu(y)
+D_{0y}^{-\beta}D_{0y}^{\alpha}\mu(y).
\end{equation}
}
{\bf Proof.}
By integrating both sides of equality (\ref{lmamud}) by $x$ 
in view of (\ref{pu3}), we obtain the condition
\begin{equation} \label{pu4}
u_x(l,y)-u_x(0,y)=D_{0y}^{\alpha}\mu(y),
\quad  0\leq y\leq T.
\end{equation}
From (\ref{lmamx0}) and (\ref{lmamxl}) by taking into account the conditions (\ref{pu2}) and (\ref{pu4})
we obtain
\begin{equation} \label{puv}
\psi(y)-2{\mathcal R}_{0y}^{0,l}\psi(\eta)=\Phi(y),
\end{equation}
with respect to the function $\psi(y)=u(0,y)+u(l,y),$  
where
$$\Phi(y)=2\sum\limits_{k=1}^{n}
\left[{\mathcal N}_{0l}^{\delta_k,0,y}+{\mathcal N}_{0l}^{\delta_k,l,y}\right]
\tau_k(\xi)+
2\left[{\mathcal R}_{0y}^{\beta,0}+{\mathcal R}_{0y}^{\beta,l}\right]
D_{0y}^{\alpha}\mu(y),$$
$\delta_k=\beta-k+1.$

Let us shou that  $y^{n-\alpha}\Phi(y)\in C[0,T].$
Since $\tau_k(x)\in C[0,l],$ then in view of the formula
(\ref{vf13}) we get
$${\mathcal N}_{0l}^{\beta-k+1,x,y}\tau_k(t)=\frac{1}{2}\int\limits_0^l
\tau_k(t)y^{\beta-k}
\phi(-\beta,\beta-k+1;-|x-t|y^{-\beta})dt\leq$$
$$\leq M \int\limits_0^l y^{\beta-k}
\phi(-\beta,\beta-k+1;-|x-t|y^{-\beta})dt=$$
$$=\frac{M}{2} y^{2\beta-k}
\phi(-\beta,2\beta-k+1;-|x-t|y^{-\beta})\Big|_{t=0}^{t=l},
$$
where $M=\max\limits_{x\in[0,l]}\tau_k(x).$
From last relation and the estimates
\begin{equation} \label{lmamofw}
\left|y^{\alpha-k} \phi(-\beta,\alpha-k+1;-xy^{-\beta})\right|
\leq C x^{-\theta}y^{\alpha-k+\beta\theta},
\end{equation}
which holds in view of (\ref{vofwr0}), follows that
$y^{k-\alpha}{\mathcal N}_{0l}^{\beta-k+1,x,y}\tau_k(\xi) \in C[0,T],$ $k=1,n.$
From $D_{0y}^{\alpha}\mu(y)\in C[0,T]$ follow that
$D_{0y}^{-\beta}D_{0y}^{\alpha}\mu(y), \,
{\mathcal R}_{0y}^{\beta,l}D_{0y}^{\alpha}\mu(y)\in C[0,T].$

The equation (\ref{puv}) is an integral Volterra equation of the second kind.
Its unique solution can be written in the form 
\begin{equation} \label{pruv}
\psi(y)=2\sum\limits_{m=0}^{\infty} {\mathcal R}_{0y}^{0,ml}\Phi(\eta).
\end{equation}
Indeed, by virtue of (\ref{ls1}) we get
$$\psi(y)-2{\mathcal R}_{0y}^{0,l}\psi(\eta)=
+2\sum\limits_{m=0}^{\infty} {\mathcal R}_{0y}^{0,ml}\Phi(\eta)
-4\sum\limits_{m=0}^{\infty} {\mathcal R}_{0y}^{0,l}{\mathcal R}_{0y}^{0,ml}\Phi(\eta)=
$$
$$=\Phi(y)+2\sum\limits_{m=1}^{\infty} {\mathcal R}_{0y}^{0,ml}\Phi(\eta)
-2\sum\limits_{m=0}^{\infty} {\mathcal R}_{0y}^{0,(m+1)l}\Phi(\eta)=$$
$$
=\Phi(y)+2\sum\limits_{m=1}^{\infty} {\mathcal R}_{0y}^{0,ml}\Phi(\eta)
-2\sum\limits_{m=1}^{\infty} {\mathcal R}_{0y}^{0,ml}\Phi(\eta)=\Phi(y).$$

We rewrite equality (\ref{pruv}) in the form 
$$\psi(y)=\Phi(y)+\int\limits_0^yW(y-\eta)\Phi(\eta)d\eta,$$
where
$$W(y)=\sum\limits_{m=1}^{\infty}
\frac{(-1)^m}{y}\phi(-\beta,0;-mly^{-\beta}).$$
From estimate (\ref{vofwr0}) follows
$$\left|y^{-1} \phi(-\beta,0;-mly^{-\beta})\right|
\leq C (ml)^{-\theta}y^{\beta\theta-1}, \quad \theta\geq -1.$$
From the last we get
$$\left|y^{1-\beta}W(y)\right|\leq \frac{C}{l^2}y^{1-\beta}y^{2\beta-1}
\sum\limits_{n=1}^{\infty}\frac{1}{m^2}=y^{\beta}\frac{C}{l^2}
\sum\limits_{n=1}^{\infty}\frac{1}{m^2}.$$
Hence, $y^{1-\beta}W(y)\in C[0,T],$ and $y^{n-\alpha}\psi(y)\in C[0,T].$

Using (\ref{pruv}) and properties  1 and 2, for the solution of equation
(\ref{puv}) we obtain the form
$$\psi(y)=
2\sum\limits_{k=1}^{n}\sum\limits_{m=0}^{\infty}
\left[{\mathcal N}_{0l}^{\delta_k,-ml,y}+{\mathcal N}_{0l}^{\delta_k,(m+1)l,y}\right]
\tau_k(\xi)+$$
$$+2\sum\limits_{m=0}^{\infty}
\left[{\mathcal R}_{0y}^{\beta,ml}+{\mathcal R}_{0y}^{\beta,(m+1)l}\right]
D_{0y}^{\alpha}\mu(y).
$$
Hence using the property 3 we obtain (\ref{ppsi}).

Now, when the function $\psi(y),$  is found, we can find $u(0,y)$ and $u(l,y),$ 
as a solution of the system
$$u(0,y)+u(l,y)=\psi(y), \quad  a_1u(0,y)+a_2u(l,y)=\varphi(y).$$
On condition $a_1\not=a_2,$ the unique solution of this system is
\begin{equation} \label{puf0}
u(0,y)=\frac{a_2}{a_2-a_1}\psi(y)+\frac{1}{a_2-a_1}\varphi(y)=\varphi_0(y),
\end{equation}
\begin{equation} \label{pufl}
u(l,y)=\frac{a_1}{a_1-a_2}\psi(y)+\frac{1}{a_1-a_2}\varphi(y)=\varphi_1(y).
\end{equation}
It is obvious that  $y^{n-\alpha}\varphi_0(y), \,y^{n-\alpha}\varphi_1(y)\in C[0,T].$

The inclusion  $y^{n-\alpha}u(x,y)\in C(\overline{\Omega})$ 
is valid if the conditions
\begin{equation} \label{puss}
\lim_{y\rightarrow 0}D_{0y}^{\alpha-n}\varphi_0(y)=\tau_n(0),
\quad
\lim_{y\rightarrow 0}D_{0y}^{\alpha-n}\varphi_1(y)=\tau_n(l).
\end{equation}
Let us show it.
From the relations (\ref{puf0}) and (\ref{pufl}) follows that
\begin{equation} \label{pduf0}
D_{0y}^{\alpha-n}\varphi_0(y)=\frac{a_2}{a_2-a_1}D_{0y}^{\alpha-n}\psi(y)
+\frac{1}{a_2-a_1}D_{0y}^{\alpha-n}\varphi(y),
\end{equation}
\begin{equation} \label{pdufl}
D_{0y}^{\alpha-n}\varphi_1(y)=\frac{a_1}{a_1-a_2}D_{0y}^{\alpha-n}\psi(y)+
\frac{1}{a_1-a_2}D_{0y}^{\alpha-n}\varphi(y).
\end{equation}
From  (\ref{ppsi}) by using (\ref{ls2}) and (\ref{ls3}) we obtain
$$D_{0y}^{\alpha-n}\psi(y)=2\sum\limits_{k=1}^{n}\sum\limits_{m=-\infty}^{\infty}
D_{0y}^{\alpha-n}
{\mathcal N}_{0l}^{\delta_k,ml,y}\tau_k(\xi)+$$
\begin{equation} \label{pdpsi}
+4\sum\limits_{m=1}^{\infty}{\mathcal R}_{0y}^{n-\beta,ml}D_{0y}^{\alpha}\mu(y)
+D_{0y}^{\beta-n}D_{0y}^{\alpha}\mu(y).
\end{equation}
The equality (\ref{pdpsi}) together with  the relation (\ref{ls5}) lead to
\begin{equation} \label{pldpsi}
\lim\limits_{y\rightarrow 0}D_{0y}^{\alpha-n}\psi(y)=
\tau_n(0)+\tau_n(l).
\end{equation}
From (\ref{pus1}), (\ref{pduf0}), (\ref{pdufl}) and (\ref{pldpsi})
follows (\ref{puss}).

Thus, a solving of Problem 1 is reduced to solving the first boundary-value problem
(\ref{lmamnu}), (\ref{puf0}), (\ref{pufl}),
for the equation (\ref{lmamud}), solution of which has the form \cite[c. 123]{PsMon}
$$u(x,y)=\sum\limits_{k=1}^{n}(-1)^{k-1}\int\limits_{0}^{l}
\tau_k(\xi)\frac{\partial^{k-1}}{\partial\eta^{k-1}}G(x,y;\xi,0)d\xi+$$
\begin{equation} \label{prpkz}
+\int\limits_{0}^{y}G_{\xi}(x,y;0,\eta)u(0,\eta)d\eta
-\int\limits_{0}^{y}G_{\xi}(x,y;l,\eta)u(l,\eta)d\eta,
\end{equation}
where
$G(x,y;\xi,\eta)=
\sum\limits_{m=-\infty}^{\infty}
\left[ \Gamma(2ml+x-\xi,y-\eta)- \Gamma(2ml-x-\xi,y-\eta)\right].$

It is obvious that (\ref{prpkz}) satisfies the equation  (\ref{lmamud}) and the conditions
(\ref{lmamnu}), (\ref{pu2}).

Let us show that (\ref{prpkz}) satisfies the condition (\ref{pu3}).
In terms of operators (\ref{lmamna}) and (\ref{lmamrn})
solution (\ref{prpkz}) can be rewritten in the form (\ref{prz1})
or
$$u(x,y)\!=\!\sum\limits_{k=1}^{n}\sum\limits_{m=1}^{\infty}
\!\left[{\mathcal N}_{0l}^{\delta_k,2ml+x,y}\!-
\!{\mathcal N}_{0l}^{\delta_k,2ml-x,y}\!+\!
{\mathcal N}_{0l}^{\delta_k,-2ml+x,y}\!-\!
{\mathcal N}_{0l}^{\delta_k,-2ml-x,y}\right]\!\tau_k(x)+$$
$$+\sum\limits_{k=1}^{n}
\left[{\mathcal N}_{0l}^{\delta_k,x,y}-{\mathcal N}_{0l}^{\delta_k,-x,y}\right]
\tau_k(x)+2\sum\limits_{m=1}^{\infty}\left[{\mathcal R}_{0y}^{0,2ml-l-x}
-{\mathcal R}_{0y}^{0,2ml-l+x}\right]\varphi_l(y)+
$$
\begin{equation} \label{pfg1}
+2\sum\limits_{m=1}^{\infty}
\left[{\mathcal R}_{0y}^{0,2ml-2l+x}-{\mathcal R}_{0y}^{0,2ml-x}\right]\varphi_0(y).
\end{equation}
We integrate the equality (\ref{pfg1}) on the interval $[0,l]$ with respect to $x.$
Using properties 5 and 6, we get
$$\int\limits_0^l u(x,y)dx=2\sum\limits_{k=1}^{n}\sum\limits_{m=0}^{\infty} (-1)^{m+1}
\left[{\mathcal N}_{0l}^{\delta_k+\beta,(m+1)l,y}+{\mathcal N}_{0l}^{\delta_k+\beta,-ml,y}\right]
\tau_k(x)+
$$
\begin{equation} \label{pfg11}
+\sum\limits_{k=1}^{n}\frac{y^{\delta+\beta-1}}{\Gamma(\delta+\beta)}
\int\limits_0^l\tau_k(x)dx+2{\mathcal R}_{0y}^{\beta,0}\psi(y)+4\sum\limits_{m=1}^{\infty}
(-1)^m{\mathcal R}_{0y}^{\beta,ml}\psi(y),
\end{equation}
where $\psi(y)=\varphi_0(y)+\varphi_l(y)$ is a solution of the Volterra equation (\ref{puv}).

From  (\ref{pruv}) we obtain
$$2{\mathcal R}_{0y}^{\beta,0}\psi(y)=
D_{0y}^{-\beta}\psi(y)=2\sum\limits_{m=0}^{\infty}
{\mathcal R}_{0y}^{\beta,ml}\Phi(y).$$
We transform the last two addents in right side of (\ref{pfg11})  with the help of
(\ref{ls1})
$$2{\mathcal R}_{0y}^{\beta,0}\psi(y)+
4\sum\limits_{m=1}^{\infty}(-1)^{m}{\mathcal R}_{0y}^{0,ml}D_{0y}^{-\beta}\psi(y)=$$
$$=2\sum\limits_{m=0}^{\infty} {\mathcal R}_{0y}^{\beta,ml}\Phi(y)+
8\sum\limits_{m=1}^{\infty}(-1)^m {\mathcal R}_{0y}^{0,ml}
\sum\limits_{s=0}^{\infty} {\mathcal R}_{0y}^{\beta,sl}\Phi(y)=$$
$$=2\sum\limits_{m=0}^{\infty} {\mathcal R}_{0y}^{\beta,ml}\Phi(y)+
8\sum\limits_{m=1}^{\infty} \sum\limits_{s=1}^{m} (-1)^s
{\mathcal R}_{0y}^{0,sl} {\mathcal R}_{0y}^{\beta,(m-s)l}\Phi(y)=$$
$$=2\sum\limits_{m=0}^{\infty} {\mathcal R}_{0y}^{\beta,ml}\Phi(y)+
4\sum\limits_{m=1}^{\infty}{\mathcal R}_{0y}^{\beta,ml}\Phi(y)
\sum\limits_{s=1}^{m} (-1)^s=$$
$$=2{\mathcal R}_{0y}^{\beta,0}\Phi(y)+
2\sum\limits_{m=1}^{\infty} {\mathcal R}_{0y}^{\beta,2ml}\Phi(y)+
2\sum\limits_{m=1}^{\infty} {\mathcal R}_{0y}^{\beta,2ml-l}\Phi(y)
-4\sum\limits_{m=1}^{\infty} {\mathcal R}_{0y}^{\beta,2ml-l}\Phi(y)=$$
\begin{equation} \label{pfg2}
=2{\mathcal R}_{0y}^{\beta,0}\Phi(y)+
2\sum\limits_{m=1}^{\infty}(-1)^m {\mathcal R}_{0y}^{\beta,ml}\Phi(y)
=2\sum\limits_{m=1}^{\infty}(-1)^m {\mathcal R}_{0y}^{\beta,ml}\Phi(y).
\end{equation}
By denoting $\mu_1(y)=D_{0y}^{-\beta}D_{0y}^{\alpha}\mu(y)$
and using the properties of operator ${\mathcal R}_{0y}^{\beta,0},$
we obtain
$$2\sum\limits_{m=1}^{\infty}(-1)^m {\mathcal R}_{0y}^{\beta,ml}\Phi(y)
=2\sum\limits_{m=1}^{\infty}(-1)^m {\mathcal R}_{0y}^{\beta,ml}\mu_1(y)
+2\sum\limits_{m=1}^{\infty}(-1)^m {\mathcal R}_{0y}^{\beta,(m+1)l}\mu_1(y)+$$
$$+2\sum\limits_{m=0}^{\infty}(-1)^m {\mathcal R}_{0y}^{\beta,ml}
2\sum\limits_{k=1}^{n}
\left[{\mathcal N}_{0l}^{\delta_k,0,y}+
{\mathcal N}_{0l}^{\delta_k,l,y}\right]\tau_k(\xi)
=D_{0y}^{-\beta}\mu_1(y)+
$$
\begin{equation} \label{pfg21}
+2\sum\limits_{k=1}^{n}\sum\limits_{m=0}^{\infty}(-1)^m
\left[{\mathcal N}_{0l}^{\delta_k+\beta,-ml,y}+
{\mathcal N}_{0l}^{\delta_k+\beta,(m+1)l,y}\right]\tau_k(\xi).
\end{equation}
By virtue of the fractional analoque of Newton-Leibniz formula \cite[p. 11]{n},
we have
$$D_{0y}^{-\beta}\mu_1(y)=D_{0y}^{-\alpha}D_{0y}^{\alpha}\mu(y)=
\mu(y)-\sum\limits_{k=1}^{n}\frac{y^{\delta_k+\beta-1}}{\Gamma(\delta_k+\beta)}
\lim\limits_{y\rightarrow 0}D_{0y}^{\alpha-k}\mu(y).$$
By taking into account last relation, 
from (\ref{pfg11}), (\ref{pfg2}) and (\ref{pfg21}) we obtain
$$\int\limits_0^l u(x,y)dx=\mu(y)+\sum\limits_{k=1}^{n}
\frac{y^{\delta_k+\beta-1}}{\Gamma(\delta_k+\beta)}
\left[\int\limits_0^l\tau_k(x)dx-
\lim\limits_{y\rightarrow 0}D_{0y}^{\alpha-k}\mu(y)\right].
$$
Thus, under matching conditions (\ref{pus2}) function $u(x,y)$ satisfies Samarskii integral condition.
The proof of Theorem 2 is complete.

\subsection{Samarskii problem for the wave equation}

Consider Problem 1 in case when  $\alpha=2,$ $a_1=1,$ $a_2=0$ and $T<l.$
In general case this problem can be solved in similar way.

{\bf Problem 2.} { Find a solution of equation 
\begin{equation} \label{pwe}
u_{xx}-u_{yy}=0,
\end{equation}
satisfying the conditions
$$u(x,0)=\tau(x), \quad u_y(x,0)=\nu(x), \quad 0<x<l,$$
$$u(0,y)=\varphi_0(y), \quad \int\limits_0^lu(x,y)dx=\mu(y), \quad 0<y<T<l.$$
where $\tau(x),$ $\nu(x),$ $\varphi_0 (y),$ $\mu (y)$ are given functions. }

Problem 2 for the wave equation was studied in the work \cite{Beylin} by the reduction to the problem
with non-local Bitsadze-Samarskii condition.
Note also the paper \cite{Gordeziani-2000} in which the non-local initial boundary value problems with
integral nonlocal boundary conditions are investigated for one-dimensional medium oscillation equations and
solutions of the corresponding problems are constructed. 
More exstensive overview of the subject of nonlocal boundary problems for wave equation 
can be found in works \cite{Beylin} and \cite {Bazhlekova-2012}.

From (\ref{lmamu11}) and (\ref{lmamu12})  by virtue of the equality
$$\int\limits_0^y[u_x(l,s)-u_x(0,s)]ds=\mu'(y)-\mu'(0),$$
we express the value of $u(l,y)$ through the data of Problem 2:
$$u(l,y)=\int\limits_0^y\nu(t)dt+\int\limits_{l-y}^l\nu(t)dt
+\tau(y)+\tau(l-y)+\mu'(y)-\mu'(0)-\varphi_0(y)\equiv \varphi_l(y).$$
Thus, Problem 2 is redused to the local first boundary value problem
for the equation (\ref{pwe}),
which solution has  the form \cite[c. 70]{Tih}
\begin{equation} \label{pswe}
u(x,y)=\frac{\tau(x+y)+\tau(x-y)}{2}+
\frac{1}{2}\int\limits_{x-y}^{x+y}\nu(t)dt+\overline{\varphi}_0(y-x)
-\overline{\varphi}_l(y+x-l),
\end{equation}
where
$\overline{\varphi}_0(y)=\varphi_0(y)H(y),$
$\overline{\varphi}_l(y)=\varphi_l(y)H(y),$
$H(y)$ is a Heaviside function, and
\begin{equation} \label{pas}
\tau(-x)=-\tau(x), \quad \tau(2l-x)=-\tau(x),\quad
\nu(-x)=-\nu(x),   \quad\nu(2l-x)=-\nu(x).
\end{equation}

Obviously, that the function (\ref{pswe}) is a solution of the equation (\ref{pwe}), 
and also that the first three conditions of Problem 2 are satisfied.

Let us show that the fourth also is satisfied.
Integrating the relation (\ref{pswe}) with respect to variable 
$x$ from $0$ to $l$
$$\int\limits_{0}^{l}u(x,y)dx=
\frac{1}{2}\int\limits_0^l[\tau(x+y)+\tau(x-y)]dx+
\frac{1}{2}\int\limits_0^l \int\limits_{x-y}^{x+y}\nu(t)dt+
\int\limits_{0}^{l}\overline{\varphi}_0(y-x)dx+$$
\begin{equation} \label{ppus}
+\int\limits_{0}^{l}\overline{\varphi}_l(y+x-l)dx
=I_1+I_2+I_3+I_4.
\end{equation}
We transform the integrals $I_k$ $(k=\overline{1,4}).$
By taking into account the equalities (\ref{pas}), we obtain
\begin{equation} \label{pi1}
2I_1=
\int\limits_y^l\tau(t)dt-\int\limits_l^{l+y}\tau(2l-t)dt
-\int\limits_{-y}^0\tau(-t)dt+\int\limits_{-y}^{l-y}\tau(t)dt=
2\int\limits_y^{l-y}\tau(t)dt,
\end{equation}
$$2I_2=
\int\limits_{-y}^{0}(t+y)\nu(t)dt+  \int\limits_{0}^{l-y}(t+y)\nu(t)dt+
\int\limits_{l-y}^{y}l\nu(t)dt+  \int\limits_{y}^{l}(l-t+y)\nu(t)dt+ $$
\begin{equation} \label{pi2}
+\int\limits_{l}^{l-y}(l-t+y)\nu(t)dt=
2\int\limits_{0}^{y}t\nu(t)dt+2\int\limits_{y}^{l-y}y\nu(t)dt
+2\int\limits_{l-y}^{l}(l-t)\nu(t)dt,
\end{equation}
\begin{equation} \label{pi3}
I_3=\int\limits_{0}^{y}\varphi_0(y-x)dx=\int\limits_{0}^{y}\varphi_0(t)dt,
\end{equation}
$$I_4=\int\limits_{l-y}^{l}\varphi_l(y+x-l)dx=\int\limits_{0}^{y}\varphi_l(t)dt
= \int\limits_{0}^{y}(y-t)\nu(t)dt
+\int\limits_{l-y}^{l}(y+t-l)\nu(t)dt+$$
\begin{equation} \label{pi4}
+\int\limits_{0}^{y}\tau(s)ds+\int\limits_{l-y}^{l}\tau(s)ds-
\int\limits_{0}^{y}\varphi_0(s)ds+\mu(y)-\mu(0)-\mu'(0)y.
\end{equation}
By virtue of (\ref{pi1}) -- (\ref{pi4}), from (\ref{ppus}) we obtain
$$\int\limits_{0}^{l}u(x,y)dx=I_1+I_2+I_3+I_4=
\left(\int\limits_{0}^{y}+\int\limits_{y}^{l-y}+\int\limits_{l-y}^{l}\right)
\tau(t)dt-\mu(0)+$$
$$+\left(\int\limits_{0}^{y}+\int\limits_{y}^{l-y}+\int\limits_{l-y}^{l}\right)
\nu(t)dt-\mu'(0)+\mu(y).$$
From the last we can see that the function
(\ref{pswe}) satisfy the condition 
$$\int\limits_{0}^{l}u(x,y)dx=\mu(y),$$
if the following matching conditions hold
$$\int\limits_{0}^{l}\tau(x)dx=\mu(0), \quad \int\limits_{0}^{l}\nu(x)dx=\mu'(0).$$
Thus the function (\ref{pswe}) is the solution of Problem 2.

\end{document}